\documentclass[12pt]{amsart}
\usepackage{times,fullpage}
\usepackage{latexsym,amscd,amssymb}
\newtheorem{thm}{Theorem}

\newtheorem{lem}[thm]{Lemma}

\theoremstyle{remark}

\theoremstyle{definition}

\newcommand{\R}{\mathbb{ R}}

\DeclareMathOperator{\DEF}{def}

\begin{document}

\title{K\"ahlerian three-manifold groups}
\author{D.~Kotschick}
\address{Mathematisches Institut, {\smaller LMU} M\"unchen,
Theresienstr.~39, 80333~M\"unchen, Germany}
\email{dieter@member.ams.org}
\date{January 7, 2013; \copyright{\ D.~Kotschick 2013}}
\thanks{Research done at the Institute for Advanced Study in Princeton with the support of The Fund For Math and The Oswald Veblen Fund.}

\begin{abstract}
We prove that if the fundamental group of an arbitrary three-manifold -- not necessarily closed, nor orientable -- is a K\"ahler group, then
it is either finite or the fundamental group of a closed orientable surface.
\end{abstract}

\maketitle



\section{Introduction}

It has been well known for more than twenty years now that the fundamental groups of compact K\"ahler manifolds are in many
ways very different from three-manifold groups. 
For example, cf.~~\cite{ABCKT}, K\"ahler groups are indecomposable under free products, and are far from real hyperbolic groups of 
dimension $\geq 3$, whereas the class of three-manifold groups is closed under taking free products and, according to Thurston,
contains many hyperbolic groups of dimension $3$. Nevertheless, it was only comparatively recently that Dimca and Suciu~\cite{DS}
proved the long-expected result that a K\"ahler group that is also the fundamental group of a closed three-manifold must be finite.
In~\cite{AIF} I gave a simple proof of that result using essentially only the Albanese map and group cohomology.

It is the purpose of this paper to give a variation of that proof that covers the fundamental groups of all three-manifolds, 
not only the closed ones. We shall prove:
\begin{thm}\label{t:main1}
If the fundamental group of some three-manifold is infinite and a K\"ahler group, then it is the fundamental group of a closed orientable
surface.
\end{thm}
This generalizes both the theorem of Dimca--Suciu~\cite{DS} and a very recent result of Friedl and Suciu~\cite{FS}, who considered
compact three-manifolds with non-empty toroidal boundary. The problem of determining all K\"ahlerian three-manifold groups was
suggested by~\cite{FS}. Discussing infinite groups only, as we do here, is an insignificant restriction, since all finite groups are in fact
K\"ahler by a classical result of Serre. Furthermore, finite three-manifold groups are well understood.

The proof of~\cite{AIF} relied on Poincar\'e duality for closed manifolds and does not work in the general case. Still,
we follow the same strategy as in that proof, looking at the homomorphism on fundamental groups induced by the Albanese map of
a K\"ahler manifold and the naturality of the cup product in group cohomology.

The proof of Theorem~\ref{t:main1} given here depends on modern developments in three-manifold topology only when dealing with 
closed manifolds with vanishing virtual first Betti numbers. The case of manifolds with non-empty boundary is more elementary, and
requires only pre-Thurston results about three-manifolds that were available forty years ago.

\section{The proof}

Let $M$ be an arbitrary three-manifold with infinite fundamental group. Assuming that $\pi_1 (M)$ is a K\"ahler group, we would 
like to prove that it is an orientable closed surface group. We do this by first going through a series of straightforward reduction steps, and then 
dealing with the crucial case of a compact aspherical three-manifold in Theorem~\ref{t:main2} below. 

\subsection{Compactness}

As usual, the term K\"ahler group denotes a group which is the fundamental group of some closed K\"ahler manifold. In particular,
K\"ahler groups are finitely presentable. Now it is a result of Jaco~\cite{Jaco} that if the fundamental group of a three-manifold is 
finitely presentable, then it is also the fundamental group of a compact three-manifold, possibly with boundary. Thus we may
assume that $M$ is compact. 

\subsection{Primeness}

We may assume that $M$ is prime, since otherwise its fundamental group would be a non-trivial free product. This is not possible, 
either by Gromov's theorem~\cite{G}, see also~\cite{ABCKT}, or, granting residual finiteness of $\pi_1(M)$, 
by~\cite[Corollary~3.2]{AIF}, where I showed that a residually finite free product has a finite index subgroup with odd first Betti number.

\subsection{Asphericity}

By assumption, $\pi_1(M)$ is infinite. It cannot be virtually cyclic since it is assumed to be K\"ahler. Therefore, our prime $M$ is 
irreducible and aspherical by the sphere theorem, compare~\cite{M}. Thus $\pi_1(M)$ is a torsion-free group of cohomological 
dimension $\leq 3$. 

\subsection{Passage to finite coverings}\label{ss:fi}

Note that finite index subgroups of K\"ahler groups are trivially K\"ahler.
Since a torsion-free group containing the fundamental group of a closed orientable surface as a finite index subgroup is itself
a closed surface group, we may replace $M$ by an arbitrary finite covering, once we check that the fundamental 
groups of non-orientable surfaces are not K\"ahler. This is indeed so, since their first Betti numbers are positive ($\R P^2$
has been excluded by the assumption that we have infinite groups) but the cup product from $H^1$ to $H^2$ is trivial, 
contradicting the Hard Lefschetz property.
Replacing $M$ by a finite covering we may assume that it is orientable, so that its boundary is orientable as well. 

\subsection{Capping off spherical boundary components}\label{ss:cap}

Next, capping off an $S^2$ in the boundary of $M$ by a three-ball does not change the fundamental group, so we may assume that
$M$ does not have spherical boundary components.

\subsection{Positivity of the first Betti number}

If $M$ has non-empty boundary, then, since the boundary is orientable and not spherical, the boundary has non-trivial first Betti number,
and so does $M$ itself by the ``half lives, half dies'' argument.

When $M$ is closed, it can of course be a rational homology sphere. By Agol's recent resolution of the virtually Haken conjecture~\cite{A},
$M$ has a finite covering with positive first Betti number. However, we do not need this recent result.
As discussed in~\cite[p.~1085/86]{AIF}, if $M$ were closed, with zero virtual first Betti number, then $\pi_1(M)$  would not be
K\"ahler, using Perelman's results and a theorem of Carlson and Toledo.

\subsection{The main argument}

We have now explained that Theorem~\ref{t:main1} follows from:
\begin{thm}\label{t:main2}
Let $M$ be a compact aspherical three-manifold with $b_1(M)>0$. If $\pi_1(M)$ is the fundamental group of a closed K\"ahler manifold $X$, 
then the Albanese map of $X$ induces an isomorphism between $\pi_1(M)$ and the fundamental group of a closed orientable surface.
\end{thm}
Of course $M$ must then have non-empty boundary. In fact, it is a classical result of Heil~\cite[Proposition~1]{Heil} that $M$ is an 
interval bundle over a surface.

For the proof of Theorem~\ref{t:main2} assume that $X$ is a closed K\"ahler manifold with $\pi_1(X)=\pi_1(M)$. The assumption 
that the first Betti number is positive implies that $X$ has a non-constant Albanese map. Since the target of the Albanese is aspherical, 
the map factors up to homotopy through the classifying space of $\pi_1(X)$, which we may take to be $M$. But the cohomological dimension 
of $M$ is strictly less than $4$, and so the Albanese image of $X$ must be a complex curve, necessarily of positive genus $g$. 
By a standard argument, this implies that the Albanese image is smooth, and the 
Albanese map has connected fibers, compare~\cite[p.~289]{F}. Therefore, we are in the situation of the following lemma:
\begin{lem}{\rm (\cite{CatF})}\label{l:C}
Let $f\colon X\longrightarrow C_g$ be a surjective holomorphic map with connected fibers from a compact
complex manifold to a curve of genus $g\geq 1$. By marking the critical values $p_1,\ldots,p_k$ of $f$ with suitable 
integral multiplicities $m_i\geq 1$, one can define the orbifold fundamental group $\pi_1^{orb}(C_g)$ of $C_g$
with respect to these multiplicities, so that one obtains a short exact sequence
\begin{equation}\label{eq:basic}
1\longrightarrow K\longrightarrow \pi_1(X)\longrightarrow \pi_1^{orb}(C_g)\longrightarrow 1
\end{equation}
in which the kernel $K$ is finitely generated, since it is a quotient of the fundamental group of a regular fiber of $f$.
\end{lem}
The point is that if there are multiplicities $m_i\geq 2$, then the orbifold fundamental group on the right is rather larger than the 
usual topological fundamental group of $C_g$, and this ensures that the kernel $K$ finitely generated, compare the discussion 
in~\cite{S,CatF,DEF}.
To prove Theorem~\ref{t:main2} we only have to prove that $K$ is trivial, for then $\pi_1(M)=\pi_1(X)= \pi_1^{orb}(C_g)$.
As $\pi_1(M)$ is torsion-free, the orbifold structure must be trivial, and $\pi_1^{orb}(C_g)$ is just the usual fundamental
group of $C_g$.

By the solution of the Fenchel conjecture, $\pi_1^{orb}(C_g)$ has a surface 
group $\pi_1(S)$ as a finite index subgroup. So at the expense of replacing $M$ by a finite 
covering, compare~\ref{ss:fi}, its fundamental group actually fits into the following extension:
\begin{equation}\label{eq:basic2}
1\longrightarrow K\longrightarrow \pi_1(M)\stackrel{\varphi}{\longrightarrow} \pi_1(S)\longrightarrow 1 \ ,
\end{equation}
with the same finitely generated $K$ as above. 
If $K=\ker (\varphi)$ is non-trivial, then by a result of Hempel and Jaco~\cite[Theorem~1]{HJ1} it is 
infinite cyclic. The five-term exact sequence of the extension~\eqref{eq:basic2} in real cohomology then reads
$$
0\longrightarrow H^1(S;\R)\stackrel{\varphi^*}{\longrightarrow} H^1(M;\R)\longrightarrow H^1(K;\R)^{\pi_1(S)}
\stackrel{\delta}{\longrightarrow} H^2(S;\R)=\R\stackrel{\varphi^*}{\longrightarrow} H^2(M;\R) \ .
$$
If the connecting homomorphism $\delta$ is non-zero, then $H^1(M;\R)=H^1(S;\R)$, with identically zero cup 
product to $H^2(M;\R)$, since the cup product is natural under $\varphi^*$, which vanishes on $H^2(S;\R)$ by 
exactness. This contradicts the K\"ahlerness of $\pi_1(M)$ via the Hard Lefschetz theorem. 
If $\delta$ is zero, then, after possibly passing to a double covering again to ensure that the action of $\pi_1(S)$ on
$H^1(K;\R)=\R$ is trivial, we have $b_1(M)=1+b_1(S)$, which is odd, again contradicting the K\"ahlerness of $\pi_1(M)$. 

This completes the proof of Theorem~\ref{t:main2} and therefore also of Theorem~\ref{t:main1}.

\subsection{An alternative argument}

The proof of Theorem~\ref{t:main2} given above has the pleasant feature of dealing with the cases that $M$ is closed or with 
non-empty boundary uniformly. In particular, it gives yet another treatment of closed three-manifolds that is different from~\cite{DS,AIF}.

Now, taking for granted the closed case, an alternative -- and much more high-tech -- treatment of manifolds with non-empty
and non-spherical boundary is implicit in my recent paper~\cite{DEF}, where I discussed K\"ahler groups of positive deficiency. 
The deficiency of a finitely presentable group is the maximum over all presentations of the difference of the number of generators 
and the number of relators. For a compact three-manifold $M$ with non-empty boundary Epstein~\cite[Lemma~2.2]{Epstein} proved that
\begin{equation*}\label{eq:Ep}
\DEF (\pi_1(M))\geq 1-\chi (M) = 1-\frac{1}{2}\chi (\partial M) \ .
\end{equation*}
Since the Euler characteristic of the boundary is non-positive, the deficiency of $\pi_1(M)$ is positive. 
If $\partial M$ has at least one boundary component with negative Euler characteristic, then the deficiency  
is at least $2$, and~\cite[Theorem~2]{DEF} applies, to say that the kernel $K$ in~\eqref{eq:basic} must be trivial and
$\pi_1(M)$ is isomorphic to $\pi_1^{orb}(C_g)$. Since $\pi_1(M)$ is torsion-free, the orbifold structure must be trivial, and this 
is an ordinary surface group.

If $M$ has toroidal boundary, then the deficiency may well be $=1$. As explained in~\cite[p.~646]{DEF}, the results 
there go through for deficiency one groups whenever one knows that the kernel $K$ in~\eqref{eq:basic} is not just finitely generated, 
but finitely presentable, or at least of type $FP_2$. This is the case here, since by a result of Scott~\cite{Scott} and 
Shalen (unpublished), $\pi_1(M)$ is coherent, meaning that any finitely generated subgroup must be finitely presentable. Thus, the 
results of~\cite{DEF} imply Theorem~\ref{t:main2} above in all cases when $M$ has non-empty boundary.

The case of a closed $M$ cannot be dealt with by appealing to~\cite{DEF}, since for a closed aspherical three-manifold
the fundamental group has vanishing deficiency by another result of Epstein, see~\cite[Lemma~3.1]{Epstein}.

\bigskip

\bibliographystyle{amsplain}

\bigskip

\end{document}